\documentclass[12pt]{amsart}
\usepackage{verbatim, amsmath,amsthm,amssymb, amscd, amsfonts}
\newtheorem{theorem}{Theorem}[section]
\newtheorem{lemma}[theorem]{Lemma}
\newtheorem{proposition}[theorem]{Proposition}
\theoremstyle{definition}
\newtheorem{definition}[theorem]{Definition}
\newtheorem{example}[theorem]{Example}
\newtheorem{rema}[theorem]{Remark}
\newtheorem{coro}[theorem]{Corollary}

\begin{document}
\newcommand{\nc}{\newcommand}
\nc{\rnc}{\renewcommand} \nc{\nt}{\newtheorem}


\nc{\TitleAuthor}[2]{\nc{\Tt}{#1}%
    \nc{\At}{#2}%
    \maketitle%
}


\nc{\Hom}{\operatorname{Hom}} \nc{\Mor}{\operatorname{Mor}} \nc{\Aut}{\operatorname{Aut}}
\nc{\Ann}{\operatorname{Ann}} \nc{\Ker}{\operatorname{Ker}} \nc{\Trace}{\operatorname{Trace}}
\nc{\Char}{\operatorname{Char}} \nc{\Mod}{\operatorname{Mod}} \nc{\End}{\operatorname{End}}
\nc{\Spec}{\operatorname{Spec}} \nc{\Span}{\operatorname{Span}} \nc{\sgn}{\operatorname{sgn}}
\nc{\Id}{\operatorname{Id}} \nc{\Com}{\operatorname{Com}}

\nc{\nd}{\mbox{$\not|$}} 
\nc{\nci}{\mbox{$\not\subseteq$}}
\nc{\scontainin}{\mbox{$\mbox{}\subseteq\hspace{-1.5ex}\raisebox{-.5ex}{$_\prime$}\hspace*{1.5ex}$}}


\nc{\R}{{\sf R\hspace*{-0.9ex}\rule{0.15ex}%
    {1.5ex}\hspace*{0.9ex}}}
\nc{\N}{{\sf N\hspace*{-1.0ex}\rule{0.15ex}%
    {1.3ex}\hspace*{1.0ex}}}
\nc{\Q}{{\sf Q\hspace*{-1.1ex}\rule{0.15ex}%
       {1.5ex}\hspace*{1.1ex}}}
\nc{\C}{{\sf C\hspace*{-0.9ex}\rule{0.15ex}%
    {1.3ex}\hspace*{0.9ex}}}
\nc{\Z}{\mbox{${\sf Z}\!\!{\sf Z}$}}


\newcommand{\gd}{\delta}
\newcommand{\sub}{\subset}
\newcommand{\cntd}{\subseteq}
\newcommand{\go}{\omega}
\newcommand{\Pa}{P_{a^\nu,1}(U)}
\newcommand{\fx}{f(x)}
\newcommand{\fy}{f(y)}
\newcommand{\gD}{\Delta}
\newcommand{\gl}{\lambda}
\newcommand{\half}{\frac{1}{2}}
\newcommand{\ga}{\alpha}
\newcommand{\gb}{\beta}
\newcommand{\gga}{\gamma}
\newcommand{\ul}{\underline}
\newcommand{\ol}{\overline}
\newcommand{\Lrraro}{\Longrightarrow}
\newcommand{\equi}{\Longleftrightarrow}
\newcommand{\gt}{\theta}
\newcommand{\op}{\oplus}
\newcommand{\Op}{\bigoplus}
\newcommand{\CR}{{\cal R}}
\newcommand{\tr}{\bigtriangleup}
\newcommand{\grr}{\omega_1}
\newcommand{\ben}{\begin{enumerate}}
\newcommand{\een}{\end{enumerate}}
\newcommand{\ndiv}{\not\mid}
\newcommand{\bab}{\bowtie}
\newcommand{\hal}{\leftharpoonup}
\newcommand{\har}{\rightharpoonup}
\newcommand{\ot}{\otimes}
\newcommand{\OT}{\bigotimes}
\newcommand{\bwe}{\bigwedge}
\newcommand{\gep}{\varepsilon}
\newcommand{\gs}{\sigma}
\newcommand{\OO}{_{(1)}}
\newcommand{\TT}{_{(2)}}
\newcommand{\FF}{_{(3)}}
\newcommand{\minus}{^{-1}}
\newcommand{\CV}{\cal V}
\newcommand{\CVs}{\cal{V}_s}
\newcommand{\slp}{U_q(sl_2)'}
\newcommand{\olp}{O_q(SL_2)'}
\newcommand{\slq}{U_q(sl_n)}
\newcommand{\olq}{O_q(SL_n)}
\newcommand{\un}{U_q(sl_n)'}
\newcommand{\on}{O_q(SL_n)'}
\newcommand{\ct}{\centerline}
\newcommand{\bs}{\bigskip}
\newcommand{\qua}{\rm quasitriangular}
\newcommand{\ms}{\medskip}
\newcommand{\noin}{\noindent}
\newcommand{\raro}{\rightarrow}
\newcommand{\alg}{{\rm Alg}}
\newcommand{\rcom}{{\cal M}^H}
\newcommand{\lcom}{\,^H{\cal M}}
\newcommand{\rmod}{\,_R{\cal M}}
\newcommand{\qtilde}{{\tilde Q^n_{\gs}}}
\nc{\e}{\overline{E}} \nc{\K}{\overline{K}} \nc{\gL}{\Lambda}
\newcommand{\tie}{\bowtie}
\newcommand {\h}{\widehat}
\newcommand {\tl}{\tilde}
\nc{\ad}{_{\dot{ad}}}
\nc{\coad}{_{\dot{{\rm coad}}}}
\nc{\ov}{\overline}
\nc{\lla}{\left\langle}  \nc{\rla}{\right\rangle}

\title{Solvability for semisimple Hopf algebras via integrals}
\author{Miriam Cohen}
\address {Department of Mathematics\\
Ben Gurion University of the Negev, Beer Sheva, Israel}
\email {mia@math.bgu.ac.il}
\author{Sara Westreich}
\address{Department of Management\\
Bar-Ilan University,  Ramat-Gan, Israel}
\email{swestric@biu.ac.il}
\thanks {This research was supported by the ISRAEL SCIENCE FOUNDATION, 170-12.}


 \begin{abstract}
We use integrals of left coideal subalgebras to develop Harmonic analysis for semisimple Hopf algebras. We show how $N^*,$ the space of functional on $N,$ is embedded in $H^*.$ We define a bilinear form on $N^*$ and show that irreducible $N$-characters are orthogonal with respect to that form. We then give an explicit formula for induced characters of $N$ and show how the induced characters are embedded in $R(H).$

In the second part we give an intrinsic definition for solvable semisimple Hopf algebras via left coideal subalgebras and their integrals. We show how this definition generalizes solvability for finite groups. In particular, commutative and nilpotent Hopf algebras ae solvable. We finally prove an analogue of Burnside theorem: A semisimple quasitriangular Hopf algebras of dimension $p^aq^b$ is solvable.
 \end{abstract}

\maketitle

\section*{Introduction}
Left coideal subalgebras of a Hopf algebra $H$ are analogues of subgroups of a group. When they  are also normal, they  give rise to Hopf quotients. Just as subgroups, left coideal subalgebras have integrals (see e.g \cite{ko,ma,sk,ta}), which when $H$ is semisimple over a field of characteristic $0$ can be chosen to be idempotent. It is this property which we use to study solvability on one hand and harmonic analysis on the other.

Generalizing harmonic analysis from finite groups to semisimple Hopf algebras, Andruskiewitsch and Natale \cite{an} worked with Hopf subalgebras playing the role of subgroups. We go a step further and establish similar results for left coideal subalgebras in the role of subgroups. This is accomplished by heavily relying on the existence of  their integrals.

Replacing normal subgroups with normal left coideal subalgebras gave a very satisfactory intrinsic definition of nilpotent Hopf algebras \cite{cwtransactions}. However, solvability remained open. The basic difficulty being that there is no obvious analogue of quotients of left coideal subalgebras over left ideals coideals.

On the level of category theory there exists a notion of solvability \cite{eno} and it is customary to define solvable Hopf algebras $H$ as those for which $Rep(H)$ is a solvable category. However, this non-intrinsic definition is unsatisfactory as it contradicts our intuition from group theory. Commutative or nilpotent Hopf algebras are not always solvable in this sense \cite[Prop.4.5(ii),Remark 4.6(i)]{eno}.

In this paper we give an intrinsic definition of solvability which is consistent with solvability for finite groups and as desired, commutative or nilpotent Hopf algebras are indeed solvable. We prove moreover  an analogue of Burnside's $p^aq^b$ theorem for semisimple quasitriangular Hopf algebras.

\medskip The paper is organized as follows. Throughout we assume $H$ is a semisimple Hpf algebra over a field $k$ of characteristic $0.$ In Section $2$ we discuss harmonic analysis for left coideal subalgebras of $H.$
Given a left coideal subalgebra $N,$ we set $B=(H^*)^N\subset H^*$ its left coideal subalgebra of invariants. Let $\gL_N$ and $\gl_B$ denote their respective integrals. We introduce a map $\gamma$ from $N^*$ to $H^*$ as follows:
$$\langle\gamma(p),h\rangle=\left\langle p,\gl_B\rightharpoonup h\right\rangle$$
for all $p\in N^*,\,h\in H.$ We discuss properties of $\gamma$ and describe in Proposition \ref{imgamma} the image $\gamma(N^*)$ inside $H^*.$

\medskip
We define a Frobenius isomorphism $F_N:\;_NN\longrightarrow\, _NN^*$
that gives rise to a symmetric form on $N^*:$
$$(p|p')_N=\lla p',F_N\minus(p)\rla,$$
which in turn yields an $N$-analogue of orthogonality for $H$-characters. We prove:

\medskip\noin{\bf Theorem \ref{thortho}}.
Let $H$ be a semisimle Hopf algebra over a field $k$ of characteristic $0$ and $N$ a  left coideal subalgebra of $H.$ Let $(\,|\,)_N$ be defined as above. Then $(\,|\,)_N$ is a non-degenerate symmetric bilinear form on $N^*$ and the irreducible characters of $N$ are orthogonal with respect to it.

\medskip As a consequence we formulate in Corollary \ref{ort} Frobenius reciprocity for characters of left coideal subalgebras in terms of $(\,|\,)_N.$ That is,
for all $\chi_i\in {\rm Irr}(H),\,\phi_j\in {\rm Irr}(N),$ we have:
$$(\chi_{i_N}|\phi_j)_N=(\chi_i,\phi_j^{\uparrow H})_H.$$

We then obtain in Proposition \ref{eq} an explicit formula for $\phi^{\uparrow H}.$
$$\phi^{\uparrow H}=\gL\coad\gamma(\phi).$$

We finally give an explicit description of the induced $N$-characters inside $R(H).$

\medskip\noin {\bf Theorem \ref{induced}}.
Let $H$ be a semisimple Hopf algebra over an algebraically closed field of characteristic $0.$ Let $N$ be a normal left coideal subalgebra of $H$ and $B=(H^*)^N,$ then
$$R(N)^{\uparrow H}=R(H)\gl_B=\{x\in R(H)|s(x)\rightharpoonup H\subset N\}.$$

In Section $3$ we define and discuss solvable Hopf algebras heavily using integrals of left coideal subalgebras. We define:

\medskip\noin{\bf Definition \ref{ds}}.
Let $H$ be a semisimple Hopf algebra. A chain of left coideal subalgebras of $H$
$$N_0\subset N_1\subset\dots\subset N_t$$
is a solvable series if for all $0\le i\le t-1,$

\medskip\noin (i)  $\gL_{N_i}\in Z(N_{i+1}),$   where $\gL_{N_i}$ is the integral of $N_i.$

\medskip\noin (ii) For all $a,b\in N_{i+1},$ $$ (a\ad b)\gL_{N_i}=\langle\gep,a\rangle b\gL_{N_i}.$$

\medskip The Hopf algebra $H$ is solvable if it has a solvable series so that $N_0=k$ and $N_t=H.$

\bigskip
We explain in  Example \ref{GG} that this definition coincides with the definition of solvability for groups. We show in  Proposition  \ref{HN} and  Proposition \ref{HK}, how natural properties of solvable groups hold also for the Hopf analogue of solvability. In particular we show:

\medskip\noin{\bf Corollary \ref{gen}}. Semisimple nilpotent Hopf algebras are solvable.

\medskip\noin This course culminates in the main theorem of this section:

\medskip\noin{\bf Theorem \ref{burnside}}. Let $H$ be a quasitriangular semisimple  Hopf algebra of dimension $p^aq^b$ over a field $k$ of characteristic $0.$  Then   $H$ is solvable. Moreover, if $N$ is a left coideal subalgebra of $H$ then $H$ has a solvable series containing $N.$

\section{Preliminaries}

Throughout this paper,  $H$ is a semisimple Hopf algebra over an algebraically closed field $k$ of characteristic $0.$ We denote by $S$ and $s$ the antipodes of $H$ and $H^*$ respectively. Let $\tl{Z}(H)$ denote the Hopf center of $H,$ that is, the maximal Hopf subalgebra of $H$ contained in  the center of $H.$

Let  $\{V_0,\dots V_{n-1}\}$  be a complete set of non-isomorphic irreducible $H$-modules. Let $\{E_0,\dots E_{n-1}\}$ and ${\rm Irr}(H)=\{\chi_0,\dots,\chi_{n-1}\}$ be the associated central primitive idempotents and  irreducible characters of $H$ respectively, where $E_0=\gL,$ the idempotent integral of $H$ and $\chi_0=\gep.$

Denote by $R(H)$  the character algebra of $H.$
Let $\dim V_i=d_i=\langle  \chi_i,1\rangle.$  Then
$$\gl=\chi_H=\sum_{i=0}^{n-1}d_i\chi_i$$
is an integral for $H^*$ satisfying $\langle\gl,\gL\rangle=1.$ By Larson \cite{la}, the irreducible characters are orthogonal with respect to the following symmetric bilinear form defined on $R(H),$
$$(\chi|\Phi)_H=\langle \chi s(\Phi),\gL\rangle$$
for all $\chi,\Phi\in R(H).$ This form can be exteded to a symmetric bilinear form on $H^*$ by:
\begin{equation}\label{larson}
(p|p')_H=\langle  s(p')p,\gL\rangle\end{equation}
for all $p,p'\in H^*.$

The Hopf algebra $H^*$ becomes a right and left $H$-module by the {\it hit} actions $\leftharpoonup$ and $\rightharpoonup$ defined for all $a\in H,\,p\in H^*,$
$$\langle  p\leftharpoonup a,a'\rangle  =\langle  p,aa'\rangle  \qquad \langle  a\rightharpoonup p,a'\rangle  =\langle  p,a'a\rangle  $$
$H$ becomes a left and right  $H^*$-module analogously.

Denote by $_{\dot{ad}} $ the left adjoint action of $H$ on itself, that is, for all $a,h\in H,$
$$h_{\dot{ad}}  a=\sum h_1aS(h_2)$$

The dual action $\coad$ of $H$ on $H^*$ is given by:
$$h\coad p=\sum S(h_2)\rightharpoonup p\leftharpoonup h_1.$$
It is straightforward to see that $\gL\coad x\in R(H)$ for all $x\in H^*$ and if $p\in R(H)$ then
\begin{equation}\label{mia}
h\coad(xp)=(h\coad x)p\end{equation}
for all $h\in H.$

Left coideal subalgebras $L$ where discussed in \cite{ko,ma,sk,ta}.  A left coideal subalgebra of $H$ is called {\it normal} if it is stable under the left adjoint action of $H.$  When $H$ is semisimple and char$(k)=0,$ then $L$ is a semisimple algebra, $H$ is free as a left and right module over $L$ and $L$ has an idempotent two sided integral $\gL_L.$ Moreover, let $L^+=L\cap\ker\gep,$ then there exists a lattice preserving isomorphism between left coideal subalgebras $L$ of $H$ and left ideals coideals $I$ of $H.$
\begin{equation}\label{recall}L \rightarrow HL^+.\end{equation}

Equivalently, $(H/HL^+)^*$ can be viewed as invariants as follows. For any subalgebra $T$ of $H^*,$ we denote by $H^{T}$ the set of $T$-invariants of $H$ under the left {\it hit} action.
That is,
\begin{equation}\label{hb}H^{T}=\{h\in H\,|b\rightharpoonup h=\langle  b,1\rangle  h,\;\forall b\in T\}
\end{equation}
It was shown in \cite[Th.2.4]{cwromania} that if $B$ is a left coideal subalgebra of $H^*$ then $N=H^B$ is a left coideal subalgebra of $H.$ In this case also $B=(H^*)^N.$

Certain normal left  coideal subalgebras appear in \cite{bu} as  a natural generalization of kernels of group representations.  For an $H$-module $V$ define its left kernel as follows:
\begin{equation}\label{lker}{\rm LKer}_V=\{h\in H\,|\,\sum h_1\ot h_2\cdot v= h\ot v,\,\forall v\in V\}.\end{equation}

For any normal left coideal subalgebra $N$ of $H$ we denote by $H//N$ the Hopf quotient $\ol{H}=H/HN^+.$ Of particular interest is the commutator algebra $H',$ defined in \cite{bu1}.  $H'$ is  the unique normal left coideal subalgebra of $H$ so that $H//H'$ is commutative and is minimal with respect to this property.  It was shown in \cite{cwcom} that $H'$ is the algebra generated by the Hopf commutators, and that
	\begin{equation}\label{tauh'}H'= H^{kG(H^*)}.\end{equation}

Summarizing  relations between algebras of invariants and their integrals, we have:

\begin{lemma}\label{basic}
Let $N$ be a left coideal subalgebra of $H$ with an idempotent integral $\gL_N.$ Then the following hold:

\medskip\noin (i) $\gL_N\rightharpoonup \gl=\gl_B\ne 0$ is an integral for $B={H^*}^N.$

\medskip\noin (ii) We have, $$\gL_N\leftharpoonup H^*=N=\gl_B\rightharpoonup H\qquad \gl_B\leftharpoonup H=B=\gL_N\rightharpoonup H^*$$

\medskip\noin (iii) $S\gL_N=\gL_N=\gl_B\rightharpoonup\gL.$ Thus $SN=H^*\rightharpoonup \gL_N.$
\end{lemma}
\begin{proof}

(i) From the definition of $(H^*)^N$ it is straightforward to see that $\gL_N\rightharpoonup H^*\subset B.$ Now, for all $b\in B,$ we have,
$$(\gL_N\rightharpoonup \gl)\langle b,1\rangle=\gL_N\rightharpoonup (\gl b)=(\gL_N\rightharpoonup \gl)b.$$
The last equality follows since $\{{\gL_N}_2\}\subset N,$ hence act trivially on $b\in B.$ Thus $\gl_B$ is a right integral for $B.$ Since $B$ is a left coideal subalgebra of $H^*$ it follows that it contains an idempotent two sided integral $\h{\gl}_B.$ It follows now that
$\lla\gl_B,1\rla\h{\gl}_B=\gl_B\h{\gl}_B=\gl_B.$ Hence $\gl_B$ is a two sided integral.

\medskip\noin (ii)    Since $N$ is a left coideal we have $\gL_N\leftharpoonup H^*\subset N.$ Now, we have $n=\frac{1}{\lla\gl_B,1\rla }\gl_B\rightharpoonup n$ for all $n\in N,$ hence $N\subset (\gl_B\rightharpoonup H).$ Finally, by part (i), $\gl_B\rightharpoonup H=\gl_B\rightharpoonup \gL\leftharpoonup H^*=\gL_N\leftharpoonup H^*.$

The second part is obtained by switching from $H$ to $H^*$ and from $N$ to $B.$

\medskip\noin (iii) By \cite{ratrace}, $\gl_B\rightharpoonup\gL=\gL_N \rightharpoonup\gl\rightharpoonup\gL=S\gL_N.$ Since $B$ is a left coidael subalgebra of $H^*,$ it follos from part (i) after switching the roles of $H$ and $H^*,$ that $S\gL_N$ is a two sided integral for $N$, hence equals $\gL_N.$

Now, by part (ii), $H^*\rightharpoonup\gL_N=H^*\rightharpoonup S\gL_N=S(\gL_N\leftharpoonup H^*)=SN.$
\end{proof}

\medskip
In what follows we present some basic results, most of them appear in the literature or follow directly after adapting suitable notations (e.g. \cite{ko,ma,sk,ta}). However, we include proofs here for the sake of completeness.

The following structure result follows from the bijective correspondence described in \eqref{recall}.
\begin{lemma}\label{1by1}
Let $N$ be a normal left coideal subalgebra of $H.$ Then there is a bijective correspondence between left coideal subalgebras of $H$ containing $N$ and left coideal subalgebras of $H//N.$
\end{lemma}
\begin{proof}
If $N\subsetneqq L,$ where $L$ is a left coideal subalgebra of $H,$ then $k\ne \pi(L)$ is a left coideal subalgebra of $H//N.$

Conversely, let $\ol{L}$ be a left coideal subalgebra of $H//N$ and let $B=(\ol{H}^*)^{\ol{L}}.$ Then $B$ is a left coideal subalgebra of $\ol{H}^*.$ and thus is a left coideal subalgebra of $H^*.$ Set $\tl{L}=H^{B},$ then $\tl{L}$ is a left coideal subalgebra of $H$ containing $N.$ By definition $B$ acts trivially on $\pi(\tl{L}),$ hence $\pi(\tl{L})\subset (\ol{H})^{B}=\ol{L}.$

We show the converse inclusion. Set $\tl{B}=(\ol{H}^*)^{\pi(\tl{L})}.$ Then $\pi(\tl{L})\supset \ol{L}$ is equivalent to $\tl{B}=(\ol{H}^*)^{\pi(\tl{L})}\subset (\ol{H}^*)^{\ol{L}}=B,$ which in turns equivalent to
$\tl{L}=H^B\subset H^{\tl{B}}.$ This result follows directly from the definition of invariants.
\end{proof}

The integrals of left coideal subalgbras contain important information about them.

\begin{lemma}\label{integral}  Let $N$ be  a left coideal subalgebra of $H$ with an  integral $\gL_N.$ Then

\medskip\noin (i) $N$ is a Hopf subalgebra if and only if $\gL_N$ is cocommutative

\medskip\noin (ii)  $N$ is a normal in $H$ if and only if $\gL_N\in Z(H).$
\end{lemma}
\begin{proof}
(i) By Lemma \ref{basic} $H^*\rightharpoonup\gL_N=S(\gL_N\leftharpoonup H^*).$ If $\gL_N$ is cocommutative then $N=SN$ hence $N$ is a right coideal as well and thus a Hopf subalgebra. The converse is trivial.

\medskip\noin (ii) By \cite[Cor. 2.5]{cwromania} $N$ is normal in $H$ if and only if $B=(H^*)^N$ is a Hopf subalgebra of $H^*.$ By part (i) this holds if and only if $\gl_B$ is cocommutative. By Lemma \ref{basic} and \cite{ratrace}, this is equivalent to $\gL_N\in Z(H).$
\end{proof}

In addition,
\begin{lemma}\label{hlamdan}
If $N$ is a left coideal subalgebra of $H$ and $B=(H^*)^N$ then the following hold:

\medskip\noin (i) $H\gL_N\cong H/HN^+$ as left $H$-modules via $\pi_{|H\gL_N},$ where $\pi$ is the natural $H$-module projection from $H$ to $H/HN^+.$

\medskip\noin (ii) $\gL_NH\cong H/N^+H$ as right $H$-modules and

\medskip\noin (iii) If $N$ is also normal in $H$ then  $\pi_{|H\gL_N}$ is an algebra isomorphism as well.
\end{lemma}
\begin{proof}
(i). We show that $\ker\pi=H(1-\gL_N).$ Clearly $H(1-\gL_N)\subset \ker\pi.$ On the other hand, if $n\in N^+,$ then $n\gL_N=0,$ and since $\gL_N$ is an idempotent, $n=n(1-\gL_N).$ It follows now that $HN^+\subset H(1-\gL_N).$

\medskip\noin (ii). The result follows since $\gL_N$ is an integral for the right coideal subalgebra $S(N)$ as well.

\medskip\noin(iii) Follows since $H/HN^+$ is an algebra, $\gL_N$ is central by Lemma \ref{integral}(ii) and thus $H\gL_N$ is an algebra as well and $\pi$ is an algebra map.
\end{proof}

\section{Harmonic analysis for left coideal subalgebras}

Motivated by \cite{an}, we introduce a  map from $N^*$ to $H^*.$
For $p\in N^*$  define $\gamma(p)\in H^*$  as follows:
\begin{equation}\label{gamma0}\left\langle\gamma(p),h\right\rangle=\left\langle p,\gl_B\rightharpoonup h\right\rangle\end{equation}
for all $h\in H.$ Since $\gl_B\rightharpoonup H\subset N$ by lemma \ref{basic}, this map is well defined. Note that the dual map $\gamma^*:H\rightharpoonup N$ is given by:
\begin{equation}\label{gamma}\gamma^*(h)=\gl_B\rightharpoonup h.\end{equation}
\begin{example}(See the discussion in \cite[(10.1)]{cr}).
Let $G$ be a finite group and let $N$ be a  subgroup of $G.$ Let $\{p_g\}$ denote the basis of $(kG)^*$ dual to the $kG$-basis $\{g\}_{g\in G}$ of $G.$  Then $B=(H^*)^N$ has a basis consisting of idempotents of the form $P_g,$ where
$$P_g=\sum_{n\in N}p_{gn}.$$
It is straightforward to check that $P_g\in (H^*)^N.$ Since $\dim B=|G|/|N|$ we have a basis. We also have,
$P_1=\sum_{n\in N}p_n$ is an integral for $B.$ So we have
$$P_1\rightharpoonup g=\left\{\begin{matrix}
                   g & {\rm if} \; g\in N \\
                   0 & {\rm otherwise}
                 \end{matrix}\right.
                 $$

Thus, for $\chi$ an $N$-character, $\gamma(\chi)$ boils down to $\tl{\chi}$ where
$$\tl{\chi}(g)=\left\{\begin{matrix}
                   \chi(g) & {\rm if} \; g\in N \\
                   0 & {\rm otherwise}
                 \end{matrix}\right.
                 $$

\end{example}

The definition of $\gamma$ implies immediately that,
\begin{coro}\label{recip1}  For all $n\in N,\;p\in N^*,$ we have $\gamma^*(n)=\lla\gl_B,1\rla  n,$   $\gamma(p)_{|N}=\lla\gl_B,1\rla  p$ and thus $\gamma(p)\rightharpoonup n =\lla\gl_B,1\rla  \sum\lla p,n_2\rla n_1.$
\end{coro}
\medskip
Recall that $N$ acts on $N^*$ by the left hit action as follows:
$$\lla n\rightharpoonup p,n'\rla=\lla p,n'n\rla $$
for all $p\in N^*,\,n,n'\in N.$

Define a left action $\star$ of $H^*$ on $N^*$ as follows:
\begin{equation}\label{cdot}\langle x\star p,n\rangle=\langle p,n\leftharpoonup x\rangle\end{equation}
for all $x\in H^*,\,p\in N^*$ and $n\in N.$ It is straightforward that this is a well defined $H^*$- module action.
\begin{lemma}\label{hstarmodule} The following hold

\medskip\noin (i) $N^*$ is a left $H^*$-module under $\star.$ The map $\gamma: {}_{H^*}N^*\rightarrow  {}_{H^*}H^*$ is a left $H^*$-module map where $H^*$ acte on itself by left multiplication.

  \medskip\noin(ii) $\gamma:_{N}N^*\rightarrow _{N}H^*$ is a left $N$-module map, where $N$ acts on $N^*$ and on $H^*$ by the left hit actions.

\medskip\noin (iii) $x\star\gep_{|N}=x_{|N}$ for all $x\in H^*.$

\medskip\noin (iv)
 $x\gl_B=\gamma(x_{|N})$ for all $x\in H^*.$ In particular, $\gl_B=\gamma(\gep_{|N}).$

 \medskip\noin (v) The map $\gamma$ is injective.
\end{lemma}
\begin{proof} (i) We have,
$$\langle x\gamma(p),h\rangle=\langle\gamma(p),h\leftharpoonup x\rangle=\langle p,\gl_B\rightharpoonup h\leftharpoonup x\rangle=\langle x\star p,\gl_B\rightharpoonup h\rangle=\langle\gamma(x\star p),h\rangle.$$

\medskip\noin (ii) For all $n\in N,\,p\in N^*,\,h\in H,$ we have,
\begin{eqnarray*}\lefteqn{\lla\gamma(n\rightharpoonup p),h\rla=}\\
&=&\lla n\rightharpoonup p,\gl_B\rightharpoonup h\rla\\
&=&\lla p,(\gl_B\rightharpoonup h)n\rla=\lla p,\gl_B\rightharpoonup (hn)\rla\qquad\text{(since $B$ acts trivially on $N$)} \\
&=& \lla \gamma(p),hn\rla\\
&=&\lla n\rightharpoonup\gamma(p),h\rla
\end{eqnarray*}

\medskip\noin (iii) For all $x\in H^*,\;n\in N,$
$$\lla x\star \gep_{|N},n\rla=\lla\gep,n\leftharpoonup x\rla=\lla x,n\rla$$

\medskip\noin (iv) Note first that for all $h\in H^*,$
$$\langle \gamma(\gep_{|N}),h\rangle=\langle \gep,\gl_B\rightharpoonup h\rangle=\langle\gl_B,h\rangle.$$
By the previous parts,
$$x\gl_B=x\gamma(\gep_{|N})=\gamma(x\star \gep_{|N})=\gamma(x_{|N})$$

\medskip\noin (v). The result follows since $\gl_B\rightharpoonup H=N.$
\end{proof}

We can describe now Im$(\gamma)$ explicitly.
\begin{proposition}\label{imgamma}
Let $H$ be a semisimple Hopf algebra over a field $k$ of characteristic $0$ and let $N$ be a left coideal subalgebra of $H.$ Let $\gamma:N^*\longrightarrow H^*$ be defined as in \eqref{gamma0}. Then
$${\rm Im}(\gamma)=H^*\gl_B
=\{x\in H^*\,|\,s(x)\rightharpoonup H\subset N\}.$$
\end{proposition}
\begin{proof}

The fact that $H^*\gl_B\subset {\rm Im}(\gamma)$ follows from Lemma \ref{hstarmodule}(iv). To see the reverse inclusion,
let $x=\gamma(p),\,p\in N^*.$ Since $\gl_B^2=\lla\gl_B,1\rla \gl_B,$ we have for all $h\in H,$
$$\langle\gamma(p)\gl_B,h\rangle=\langle \gamma(p),\gl_B\rightharpoonup h\rangle=\lla\gl_B,1\rla \langle p,\gl_B\rightharpoonup h\rangle=\lla\gl_B,1\rla \langle \gamma(p), h\rangle=\lla\gl_B,1\rla \langle x, h\rangle.$$
Hence $\langle\gl_B,1\rangle x=\gamma(p)\gl_B\in H^*\gl_B.$ Thus $x\in H^*\gl_B.$

We show the second equality. By Lemma \ref{basic}(ii) and (iii) we have for all $x=y\gl_B\in H^*\gl_B,$
\begin{equation}\label{y}s(x)\rightharpoonup H=s(y\gl_B)\rightharpoonup H=\gl_Bs(y)\rightharpoonup H\subset N.\end{equation}
Conversely, assume $s(x)\rightharpoonup H\subset N,$ then in particular $s(x)\rightharpoonup \gL\in N.$  It follows that $$\gl_Bs(x)\rightharpoonup\gL=\lla\gl_B,1\rla s(x)\rightharpoonup\gL.$$ Thus $\gl_Bs(x)=\lla\gl_B,1\rla s(x).$ Hence we obtain
\begin{equation}\label{x}x=\frac{1}{\lla\gl_B,1\rla }x\gl_B\in H^*\gl_B.\end{equation}
\end{proof}

\begin{rema}
 Lemma \ref{hstarmodule} and the first equality in Proposition \ref{imgamma}, imply that $N^*\cong H^*\gl_B$ as left $H^*$-modules. This is in fact a dual version of Lemma \ref{hlamdan}.
\end{rema}

\medskip
As a semisimple algebra, let  $\{W_0,\dots W_{m-1}\}$  be a complete set of non-isomorphic irreducible $N$-modules. Let $\{T_0,\dots T_{m-1}\}$ and ${\rm Irr}(N)=\{\phi_0,\dots,\phi_{m-1}\}$ be the associated central primitive idempotents and  irreducible characters of $N$ respectively, where $T_0=\gL_N,$ the idempotent integral of $N$ and $\phi_0=\gep_{|N}.$

 Let $\{t_j\}$ be the set of primitive idempotents in $N$ so that $T_it_j=\gd_{ij}.$ Then, $W_j\cong Nt_j.$
 Since as $N$-modules we have $NT_j=(Nt_j)^{\oplus\lla\phi_j,1\rla},$ it follows that for all $n\in N,$
 \begin{equation}\label{phit}\lla\phi_j,nT_i\rla=\gd_{ij}\lla\phi_j,n\rla.\end{equation}
 Note that.
 \begin{equation}\label{red}
 {\chi_i}_{|N}=\sum_j m(\phi_j,{\chi_i}_{|N})\phi_j=\sum_j \left\langle\chi_i,t_j\right\rangle\phi_j,
 \end{equation}
 where $m(\phi_j,{\chi_i}_{|N})$ is the number of appearances of $\phi_j$ in ${\chi_i}_{|N}.$ Thus  $\lla\chi_i,t_j\rla$ are non-negative integers. Indeed, since $N$ is a semisimple algebra we have $\chi_{i_{|N}}=\sum_j m(\phi_j,{\chi_i}_{|N})\phi_j.$  Applying both sides to $t_j$  yields
$m(\phi_j,{\chi_i}_{|N})=\lla \chi_i,t_j\rla.$
which implies \eqref{red}.

As usual, for any $N$-module $V$ with character $\phi,$ we denote by $V^H$ the induced $H$-module $H\ot_N V$ and by $\phi^{\uparrow H}$ its induced character.

Denote by $R(N)$ the $k$-span of characters of irreducible representations of $N.$ Then \eqref{phit} implies that $Z(N)$ and $R(N)$ are dual spaces with dual bases $\{\frac{1}{\lla\phi_j,1\rla}T_j\}$ and $\{\phi_j\}.$

Extending linearly the induction of characters we define  $\phi^{\uparrow H}=\sum\ga_j\phi_j^{\uparrow H}$ for all $\phi=\sum\ga_j\phi_j\in R(N),$ and $R(N)^{\uparrow H}=\{\phi^{\uparrow H},\,\phi\in R(N)\}.$

 Variations of Frobenious reciprocity for algebras are used in the literature. We use Radford's trace formula for Hopf algebras \cite{ratrace} to show the following variation.
\begin{proposition}\label{starting}
Let $N$ be a left coideal subalgebra of $H.$ Then
 $$\phi_j^{\uparrow H}=\sum_i \left\langle\chi_i,t_j\right\rangle \chi_i.$$
  In particular,  $m(\chi_i,\phi_j^{\uparrow H}),$ the number of  appearances  of $\chi_i$ in $\phi_j^{\uparrow H},$ equals the number of appearances of $\phi_j$ in ${\chi_i}_{|N},$ and moreover, this number equals $\lla\chi_i,t_j\rla.$ That is,
  $$\lla\chi_i,t_j\rla=m(\chi_i,\phi_j^{\uparrow H})=m(\phi_j,{\chi_i}_{|N}).$$
\end{proposition}
\begin{proof}
Since $H$ is free as a right $N$-module we have $$V_{\phi_j}^{\uparrow H}=H\ot_N Nt_j\cong Ht_j$$
as left $H$-modules. By the trace formula \cite{ratrace} we have for all $h\in H,$
$$ \left\langle\phi_j^{\uparrow H},h\right\rangle =\sum \left\langle \gl,h\gL_1t_jS\gL_2\right\rangle =\sum_i\left\langle\chi_i,t_j\right\rangle \left\langle\chi_i,h\right\rangle$$
The last equality follows since we have $\gL\ad t_j\in Z(H,\; \gl=\sum d_i\chi_i$ and $\chi_i\leftharpoonup z=\frac{1}{d_i}\left\langle\chi_i,z\right\rangle\chi_i$ for all $z\in Z(H)$ (see e.g.\cite{cwcom}).

The desired result follows now from \eqref{red}.
\end{proof}

As a corollary we obtain,
\begin{coro}\label{corstarting}
Let $N$ be a left coideal subalgebra of $H.$ Then the following hold:

\medskip\noin (i) $(\dim B)\lla\phi_j,1\rla=\sum_i d_i\left\langle\chi_i,t_j\right\rangle.$

\medskip\noin (ii) $\gl_{|N}=(\dim B)\sum_j\lla\phi_j,1\rla\phi_j.$

\medskip\noin (iii) If $\gL_N$ is an idempotent integral of $N$ and $\gl_B=\gL_N\rightharpoonup\gL,$ then $\lla\gl_B,1\rla  =\dim B.$
\end{coro}
\begin{proof}
\medskip\noin (i) The result follows by applying both sides of Proposition \ref{starting}  to $1$ and since $\dim (H\ot_N Nt_j)=(\dim B)\lla\phi_j,1\rla.$

\medskip\noin (ii)   The result follows from  \eqref{red} and (i) applied to $\gl=\sum d_i\chi_i.$

\medskip\noin (iii) By Lemma \ref{basic}(i), $\lla\gl_B,1\rla =\lla\gl,\gL_N\rla.$ By part (ii), this equals in turn   to
$(\dim B)\sum_j\lla\phi_j,1\rla\lla\phi_j,\gL_N\rla.$ The fact that $\gL_N=T_0,$ a primitive central idempotent of $N$ with $\phi_0=\gep_{|N},$ yields $\lla\gl_B,1\rla=\dim B.$
\end{proof}

\medskip
We define a Frobenius map $F_N:\;_NN\longrightarrow\, _NN^*$ as follows:
$$F_N(n)=\frac{1}{\lla\gl_B,1\rla }n\rightharpoonup\gl_{|N}.$$
We show,
\begin{lemma}\label{inv}

\noin (i) The map $F_N$ is self adjoint and  bijective.
Its inverse  map $F_N\minus:N^*\longrightarrow N$ is given by:
$$F_N\minus (p)=\frac{1}{\lla\gl_B,1\rla }S\left(\gamma(p)\rightharpoonup \gL_N\right)=\frac{1}{\lla\gl_B,1\rla }\gL_N\leftharpoonup s(\gamma(p))$$
for all $p\in N^*.$

\medskip\noin (ii) For a central primitive idempotent of $N,\,T_k,$ we have:
$$F_N(T_k)=\lla\phi_k,1\rla\phi_k.$$
\end{lemma}
\begin{proof}  (i) $F_N=F_N^*$ since $\gl$ is cocommutative. Now, since $S\gL_N=\gL_N$ and $s(\gl_B)=\gl_B,$ we have for all $n\in N,$
\begin{eqnarray*}\lefteqn{F_N\minus F_N(n)=}\\
&=&\frac{1}{\lla\gl_B,1\rla ^2}S\left(\gamma(n\rightharpoonup\gl_{|N})\rightharpoonup \gL_N\right)\\
&=& \frac{1}{\lla\gl_B,1\rla ^2}S\left((n\rightharpoonup\gamma(\gl_{|N}))\rightharpoonup \gL_N\right)\quad(\text{by Lemma \ref{hstarmodule}(ii)})\\
&=& \frac{1}{\lla\gl_B,1\rla ^2}S\left((n\rightharpoonup\gl\gl_B)\rightharpoonup \gL_N\right)\qquad(\text{by Lemma \ref {hstarmodule}(iv)})\\
&=& \frac{1}{\lla\gl_B,1\rla }S\left((n\rightharpoonup\gl)\rightharpoonup \gL_N\right)\\
&=& \frac{1}{\lla\gl_B,1\rla }S\left((n\rightharpoonup\gl)\rightharpoonup\gL\leftharpoonup \gl_B\right)\qquad(\text{by lemma \ref{basic}(iii)})\\
&=&\frac{1}{\lla\gl_B,1\rla }S\left( Sn\leftharpoonup \gl_B\right)\qquad\qquad\qquad(\text{by \cite{ratrace}})
\\
&=&n.
\end{eqnarray*}
Hence $F_N$ is injective and thus bijective.

\medskip\noin(ii) For all $n\in N,$
\begin{eqnarray*}\lefteqn{\lla F_N(T_k),n\rla=}\\
&=&\frac{1}{\lla\gl_B,1\rla }\lla T_k\rightharpoonup \gl_{|N},n\rla=\\
&=&\frac{1}{\lla\gl_B,1\rla }\lla  \gl,nT_k\rla\\
&=&\sum_j\lla\phi_j,1\rla\lla\phi_j,nT_k\rla\qquad(\text{By Lemma \ref{starting}(iii)})\\
&=&\lla\phi_k,1\rla\lla\phi_k,n\rla\qquad(\text{by \eqref{phit}}).
\end{eqnarray*}
Hence the result follows.
\end{proof}

The map $F_N$ gives rise to a symmetric form on $N^*$ by:
\begin{equation}\label{form}(p|p')_N=\lla p',F_N\minus(p)\rla=\frac{1}{\lla\gl_B,1\rla }\lla s\gamma(p)\star p',\gL_N\rla,\end{equation}
where $(\star)$ is the action defined in \eqref{cdot}.
\begin{rema}\label{subalgebra}
If $N$ is a Hopf subalgebra of $H,$ then  the symmetric form defined in \eqref{form} boils down to the bilinear form on $H^*$ given in \eqref{larson}.

Indeed, when $N$ is a Hopf subalgebra then $N^*$ is a Hopf algebra as well, and we obtain,
$$\lla p',F_N\minus(p)\rla=\frac{1}{\lla\gl_B,1\rla }\langle p',S(\gamma(p)\rightharpoonup\gL_N\rangle=\langle p',S(p\rightharpoonup\gL_N)\rangle=\langle p's(p),\gL_N\rangle.$$
The second equality follows from Corollary \ref{recip1}.
\end{rema}
As for Hopf algebras we obtain orthogonality of characters of left coideal subalgebras.
\begin{theorem}\label{thortho}
Let $H$ be a semisimle Hopf algebra over a field $k$ of characteristic $0$ and $N$ a  left coideal subalgebra of $H.$ Let $(\,|\,)_N$ be defined as in \eqref{form}. Then $(\,|\,)_N$ is a non-degenerate symmetric bilinear form on $N^*$ and the irreducible characters of $N$ are orthogonal with respect to it.
\end{theorem}
\begin{proof}
The form is symmetric since $F_N$ (and thus $F_N\minus$) is self adjoint. It is non-degenerate since $F_N$ is bijective. Now, by Lemma \ref{inv}(ii) and \eqref{phit},
$(\phi_j|\phi_l)_N=\lla \phi_l, \frac{1}{\lla\phi_j,1\rla}T_j\rla=\gd_{lj}.$
\end{proof}

 As a consequence we formulate a Frobenius reciprocity  for characters of left coideal subalgebras in terms of $(\,|\,)_N.$

\begin{coro}\label{ort}
For all $\chi_i\in {\rm Irr}(H),\,\phi_j\in {\rm Irr}(N),$ we have:
$$(\chi_{i_N}|\phi_j)_N=(\chi_i|\phi_j^{\uparrow H})_H=\lla\chi_i,t_j\rla.$$

Consequently,
$$m(\phi_j|{\chi_i}_{|N})=(\chi_{i_N}|\phi_j)_N =(\chi_i|\phi_j^{\uparrow H})_H=m(\chi_i|\phi_j^{\uparrow H}).$$
\end{coro}
\begin{proof}
 By Lemma \ref{starting}(i) and (iii) and  orthogonality of irreducible characters of $H$ (and of $N$)  with respect to $(\,,\,)_H,$ (respectively $(\,|\,)_N$), each of the first two expressions equals $\lla\chi_i,t_j\rla.$

 The second formula follows from Proposition \ref{starting}.
\end{proof}
We obtain also an explicit description of $\phi^{\uparrow H}.$
\begin{proposition}\label{eq}
Let $N$ be a  left coideal subalgebra of $H$ and let  $\phi\in R(N).$ Then $$\gL\coad\gamma(\phi)=\phi^{\uparrow H}.$$
\end{proposition}
\begin{proof}
let $\phi_j$ be an irreducible character of $N.$ We show that for all $\chi_i\in R(H),$
$$(\chi_i,\phi_j^{\uparrow H})_H=(\chi_i,\gL\coad\gamma(\phi_j))_H.$$ Indeed,
\begin{eqnarray*}\lefteqn{\lla\phi_j,1\rla(\chi_i,\gL\coad\gamma(\phi_j))_H=}\\
&=& \lla \gL\coad\gamma(F_N(T_j))s(\chi_i),\gL\rla\qquad\qquad\quad(\text{by Lemma \ref{inv}(ii)})\\
&=& \lla \gL\coad(\gamma(F_N(T_j))s(\chi_i)),\gL\rla\qquad\qquad(\text{by \eqref{mia}})\\
&=& \lla F_N(T_j),\gl_Bs(\chi_i)\rightharpoonup\gL\rla\qquad\qquad\qquad(\text{since $\gL$ is cocommutative})\\
&=& \lla T_j,F_N(\gl_Bs(\chi_i)\rightharpoonup\gL)\rla\qquad\qquad\qquad\;(\text{since $F_N$ is self adjoint})\\
&=& \lla\gl_B,1\rla \minus\lla T_j,(\gl_Bs(\chi_i)\rightharpoonup\gL)\rightharpoonup\gl\rla\;(\text{by definition of $F_N$ and since $T_j\in N$})\\
&=& \lla\gl_B,1\rla \minus\lla T_j,\chi_i\gl_B\rla\qquad\qquad\qquad\quad(\text{by \cite{ratrace}})\\
&=&\lla T_j,\chi_i\rla\qquad\qquad\qquad\qquad\qquad\qquad\;(\text{since }T_j\in N)\\
&=&\lla\phi_j,1\rla\lla t_j,\chi_i\rla\qquad\qquad\qquad\qquad\qquad(\text{by \eqref{phit} and \eqref{red}})\\
&=&\lla\phi_j,1\rla(\chi_i,\phi_j^{\uparrow H})_H\qquad\qquad\quad\qquad\qquad\;(\text{by Lemma \ref{ort}})
\end{eqnarray*}
Since the bilinear form $(\,,\,)_H$ is non-degenerate on $R(H),$ the result follows.
\end{proof}
In particular we obtain,
\begin{coro}\label{lamdanh}
$(\gep_{|N})^{\uparrow H}=\gL\coad \gl_B.$
\end{coro}
More can be said if $N$ is also normal in $H.$
\begin{coro}\label{co} Let $N$ be a normal left coideal subalgebra of $H$ and $\Psi\in R(H).$ Then $(\Psi_{|N})^{\uparrow H}=\Psi\gl_B.$ In particular $(\gep_{|N})^{\uparrow H}=\gl_B.$
\end{coro}
\begin{proof}
By Lemma \ref{hstarmodule}(iv), we have $\gamma(\Psi_{|N})=\Psi\gl_B.$  Hence by Proposition \ref{eq} we have that $(\Psi_{|N})^{\uparrow H}=\gL\coad (\Psi\gl_B).$ Since $B$ is a Hopf subalgebra of $H$ it follows that $\gl_B$ is cocommutative, hence so is $\Psi\gl_B$ and the result follows.
\end{proof}
(The corollary above could also be deduced from \cite[Cor. 3.2, 3.11]{bu1}).

\medskip
Based on Theorem \ref{imgamma} we obtain a characterization of  the image of $R(N)^{\uparrow H}$ inside $R(H).$
\begin{theorem}\label{induced}
Let $H$ be a semisimple Hopf algebra over an algebraically closed field of characteristic $0.$ Let $N$ be a normal left coideal subalgebra of $H$ and $B=(H^*)^N,$ then
$$R(N)^{\uparrow H}=R(H)\gl_B=\{x\in R(H)|s(x)\rightharpoonup H\subset N\}.$$
\end{theorem}
\begin{proof}
We show first that $R(N)^{\uparrow H}=\{\phi^{\uparrow H}\,|\,\phi\in R(N)\}\subset R(H)\gl_B.$  By Proposition \ref{eq} and Theorem \ref{imgamma}, we have,
$$\phi^{\uparrow H}=\gL\coad\gamma(\phi)=\gL\coad(y\gl_B)=(\gL\coad y)\gl_B\in R(H)\gl_B,$$
where $y\in H^*.$

The reverse inclusion is Corollary \ref{co}.
The fact that $R(H)\gl_B\subset \{x\in R(H)|s(x)\rightharpoonup H\subset N\}$ is just a particular case of \eqref{y}, where $x=\Psi\gl_B\in R(H).$

We show now the reverse inclusion. Assume $x\in R(H)$ satisfies $s(x)\rightharpoonup H\subset N.$ By  \eqref{x} we have $$x=\frac{1}{\lla\gl_B,1\rla }x\gl_B\in R(H)\gl_B.$$
\end{proof}

\section{Solvablity for Hopf algebras}

In this section we suggest an intrinsic definition for solvable Hopf algebra. We then justify it based on its properties.

Observe first that
if $L$ and $N$ are left coideal subalgebras of $H,\;B_L=(H^*)^L,\;B_N=(H^*)^N$ are the corresponding left coideal subalgebras of $H^*,$  then the algebra $B_LB_N$ generated by $B_L\cup B_N$ is a left coideal subalgebra of $H^*$ and
\begin{equation}\label{cap}(H^*)^{L\cap N}=B_LB_N.\end{equation}

Indeed, since $L\cap N\subset L$ (respectively $N$), it follows that $B_L\subset (H^*)^{L\cap N}.$ Also  $B_N\subset (H^*)^{L\cap N}.$ Since $(H^*)^{L\cap N}$ is an algebra it must contain the algebra generated by $B_L\cup B_N,$ that is $B_LB_N\subset (H^*)^{L\cap N}.$ Conversely,   $B_L\subset B_LB_N$ implies that $H^{B_LB_N}\subset H^{B_L}=L.$ The same for $N,$ and so $H^{B_LB_N}\subset L\cap N$ implying that $(H^*)^{L\cap N}\subset B_LB_N.$
We show now,
\begin{lemma}\label{integrals}
Let $N$ and $L$ be  left coideal subalgebra of $H$ and let $\gl_{B_N},\,\gl_{B_L}$ denote the respective non-zero integrals of the left coideal subalgebras $B_N$ and $B_L.$

\medskip\noin (i) If $\gl_{B_N}\gl_{B_L}=\gl_{B_L}\gl_{B_N}$ then $\gl_{B_L}\gl_{B_N}$ is an integral for $B_LB_N.$

\medskip\noin (ii) The condition in (i) is satisfied if either $N$ is a normal Hopf subalgebra of $H,$  or $R(H)$ is commutative and $N$ and $L$ are normal left coideal subalgebras of $H.$
\end{lemma}
\begin{proof}
\medskip\noin (i) If $\gl_{B_N}\gl_{B_L}=\gl_{B_L}\gl_{B_N}$ then $\gl_{B_L}\gl_{B_N}b=\langle b,1\rangle \gl_{B_L}\gl_{B_N}$ for all $b\in B_N\cup B_L,$ and thus for all $b'$ in the algebra generated by $B_N\cup B_L.$ Since $\langle \gl_{B_N}\gl_{B_L},1\rangle\ne 0,$ uniqueness of the integral implies the result.

\medskip\noin (ii) If $L$ is a normal Hopf subalgebra of $H$ then $B_L$ is a normal Hopf subalgebra of $H^*$ and by Lemma \ref{integral} $\gl_{B_L}\in Z(H^*).$ If $N$ and $L$ are normal left coideal subalgebras then  $\gl_{B_N},\,\gl_{B_L}\in R(H),$ and commutativity of $R(H)$  implies that $\gl_{B_N}\gl_{B_L}=\gl_{B_L}\gl_{B_N}.$ In both cases (i) is satified.
\end{proof}

For groups, let $K$ be a normal subgroup of $G$ and $\pi:G \rightarrow G/K.$ If $L$ is a subgroup of $G$ so that $L\cap K=1$ then $\pi_{|L}$ is injective.  Surprisingly, this idea is not necessarily true for Hopf algebras if we replace $K$ by a normal left coideal subalgebra $N.$ However, it does hold in the following:
\begin{proposition}\label{embedding0} Let $N$ be a normal left coideal subalgebra of $H,\;\pi:H\rightarrow H//N$ and  $L$ any  left coideal subalgebra of $H.$  Assume $\gl_{B_N}\gl_{B_L}=\gl_{B_L}\gl_{B_N}.$ If $L\cap N=k$ then $\pi_{|L}$ is injective
\end{proposition}
\begin{proof}
 If $L\cap N=k$ then by \eqref{cap} $B_LB_M=(H^*)^k=H^*.$ By assumption and Lemma \ref{integrals} $\gl_{B_L}\gl_{B_N}=\gl.$

Now, let $a\in L$ and assume $\gL_Na=0.$ We have by Lemma \ref{basic}(i), $\gl_{B_N}\rightharpoonup\gL=\gL_N.$ Hence,
$$0=\gl_{B_L}\rightharpoonup(\gL_Na)=(\gl_{B_L}\rightharpoonup\gL_N)a=(\gl_{B_L}\gl_{B_N}\rightharpoonup\gL)a= (\gl\rightharpoonup\gL)a=a.$$
The second equality follows from the fact that $B_L$ is a left coideal subalgebra acting trivially on $L.$ By Lemma \ref{integral}(ii),  $\gL_N$ is central in $H,$ hence we obtain $a\gL_N=0$ implies $a=0$ for all $a\in L.$ The desired result follows now from Lemma \ref{hlamdan}.
\end{proof}
In particular, by Lemma \ref{integrals}(ii),
\begin{coro}\label{embedding}
Let $N$ be a normal left coideal subalgebra of $H$ and $\pi:H\rightarrow H//N.$ Assume $L\cap N=k.$ If either $N$ is a normal Hopf subalgebra of $H,$  or $R(H)$ is commutative and  both $N$ and $L$ are normal left coideal subalgebras, then $\pi_{|L}$ is injective.
\end{coro}
The assumptions in Lemma \ref{integrals} and in Corollary \ref{embedding} are necessary. We shall see that in the following counter example due to Skryabin, whom we wish to thank for this contribution.
\begin{example} Let $H=(kS_3)^*,\;B_N=k(1,2),\;B_L=k(13)\subset H^*.$
Note first that $B_NB_L=H^*.$  Clearly $R(H)=kS_3$ is not commutative.
Let $N=H^{B_N},\,L=H^{B_L},$ then $\dim N=\dim L=3.$

Since $B_N$ is not normal in $H^*,$ we know that $N$ is not a Hopf subalgebra of $H.$ Since $B_NB_L=H^*,$ it follows from \eqref{cap} that $N\cap L=k.$ On the other hand, since $B_N\cong (H/HN^+)^*,$ it follows that $\dim HN^+=4.$ Thus $(\dim L+\dim HN^+)>6,$ implying that $L\cap HN^+\ne \{0\}.$  This shows that the restriction of $\pi:H\rightarrow H/HN^+$ to $L$   is not injective.   Thus the conditions in Corollary \ref{embedding} are necessary.

Indeed, this example violates the condition of Lemma \ref{integrals} and thus of Proposition \ref{embedding0}, as
$$\gl_{B_N}\gl_{B_L}=1+(12)+(13)+(132)\ne\gl_{B_L}\gl_{B_N}=1+(12)+(13)+(123)$$
None of the products is an integral for $B_LB_M=H^*=kS_3.$
\end{example}

\bigskip
We are ready to define solvability for Hopf algebras.

\begin{definition}\label{ds}
Let $H$ be a semisimple Hopf algebra. A chain of left coideal subalgebras of $H$
$$N_0\subset N_1\subset\dots\subset N_t$$
is a solvable series if for all $0\le i\le t-1,$

\medskip\noin (i)  $\gL_{N_i}\in Z(N_{i+1}),$   where $\gL_{N_i}$ is the integral of $N_i.$

\medskip\noin (ii) For all $a,b\in N_{i+1},$ $$ (a\ad b)\gL_{N_i}=\langle\gep,a\rangle b\gL_{N_i}.$$

\medskip The Hopf algebra $H$ is solvable if it has a solvable series so that $N_0=k$ and $N_t=H.$
\end{definition}

\begin{example}\label{GG}
For finite groups, the definition of solvable Hopf algebras boils down to the usual definition of solvability for groups.

To see this, let $G'\subset G$ and $\gL_{G'}=\frac{1}{|G'|}\sum_{x\in G'}x.$
 If $\gL_{G'}$ is central in $G,$ then for $g\in G$ we have,
 $$\frac{1}{|G'|}\sum_{x\in G'}xg=\gL_{G'}g=g\gL_{G'}=\frac{1}{|G'|}\sum_{x\in G'}gx.$$
 It follows that for $x\in G',\;xg=gy$ for some $y\in G'.$ Thus $G'$ is normal in $G.$

We claim now that $G/G'$ is abelian if and only if $(a\ad b)\gL_{G'}=b\gL_{G'},$ for all $a,b\in G.$ Indeed, assume $\ol{a}\ol{b}=\ol{b}\ol{a}$ for all $\ol{a},\ol{b}\in G/G'.$ Then $aba\minus= by$ for some $y\in G'.$ This implies that $$aba\minus \sum_{x\in G'}x =by\sum_{x\in G'}x=b\sum_{x\in G'}x.$$

Conversely, if $aba\minus\sum_{x\in G'}x=b \sum_{x\in G'}x$ then $aba\minus=by$ for some $y\in G'.$ This shows our claim.
\end{example}

\medskip By definition we obtain that homomorphic images of solvable Hopf algebras are solvable. Moreover,
\begin{proposition}\label{HN}
Let $N$ be a normal left coideal subalgebra of $H.$ If the Hopf algebra $H//N$ is solvable, then $H$ has a solvable series with $N_0=N,\;N_t=H.$
\end{proposition}
\begin{proof}
Set  $\ol{H}=H//N=\pi(H)$ and let  $k\subset \ol{N_1}\subset\dots\subset\ol{N_t}=\ol{H}$ be a solvable series for $\ol{H}.$ By using Lemma \ref{1by1}, let $N_i$ be the left unique coideal subalgebra of $H$ containing $N$ so that $\pi(N_i)=\ol{N}_i.$ We want to show that $N\subset N_1\subset\dots N_t=H$ is a solvable series in $H.$

Let $\gL_{N_i}$ denote the integral of $N_i,$ then $\pi(\gL_{N_i})$ is an integral for $\ol{N}_i,$  which by assumption is central in $\ol{N}_{i+1}.$ Hence,
$$\pi(a\gL_{N_i})=\pi(a)\pi(\gL_{N_i})=\pi(\gL_{N_i})\pi(a)=\pi(\gL_{N_i}a)$$
for all $a\in N_{i+1}.$
  Since $N\subset N_i,$ we have $\gL_{N_i}\gL_{N}=\gL_{N_i},$ and since $N$ is normal we have $\gL_N\in Z(H).$ Therefore,
$$\pi(a\gL_{N_i}\gL_N)=\pi(\gL_{N_i}a\gL_N).$$
 By Lemma \ref{hlamdan} $\pi$ is injective on $H\gL_N,$ hence we have
$$a\gL_{N_i}=\gL_{N_i}a$$
for all $a\in N_{i+1}.$
This proves (i) in the definition of solvability.

Similarly, by assumption we have for $a,b\in N_{i+1},$
$$\pi\left((a\ad b)\gL_{N_i}\right)=\left(\pi(a)\ad\pi(b)\right)\pi(\gL_{N_i})=\langle \gep,a\rangle \pi(b\gL_{N_i}).$$
As before, this implies that
$$\pi\left((a\ad b)\gL_{N_i}\gL_N\right)=\pi((a\ad b)\gL_{N_i})\pi(\gL_N)=\langle \gep,a\rangle \pi(b\gL_{N_i}\gL_N),$$
and by Lemma \ref{hlamdan},
$(a\ad b)\gL_{N_i}=\langle \gep,a\rangle b\gL_{N_i}.$ Thus (ii) is satisfied.
\end{proof}
Commutative Hopf algebras are solvable by definition with $k\subset H$ as a solvable series. We wish to show that, as for groups, nilpotent Hopf algebras are solvable as well.

Recall,
\cite[Definition 1.2 (7) and (8), Lemma 1.3]{cwtransactions},  an ascending central series for $H$ is a series of normal left coideal subalgebras
\begin{equation}\label{nil}k\subset Z_1\subset  \cdots\subset Z_t\subset\cdots \end{equation}
so that  $Z_1=\tl{Z}(H)$ is the Hopf center of $H,\;H_0=k,\;H_1=H//\tl{Z}(H),$   and   $Z_i=H^{co\pi_i},$ where
$$\pi_{i}:H\rightarrow H_{i-1}//\tl{Z}(H_{i-1})=H_{i}$$ for all $i\ge 0.$ Note that we have also $H_i=H//Z_i.$

$H$ is {\it nilpotent} if $H_t=H$ for some $t.$

Based on our approach here, we suggest below another characterization for nilpotent Hopf algebras.
\begin{proposition}\label{nilnew}
A semisimple Hopf algebra $H$ is nilpotent if and only if it has a series of normal left coideal subalgebras
$$k=N_0\subset N_1\subset\cdots\subset N_t=H$$ so that $$N_{i+1}\gL_{N_i}\subset Z(H\gL_{N_i}),$$
for all $1\le i\le t,$ where $Z(H\gL_{N_i})$ is the  center of $H\gL_{N_i}.$
\end{proposition}
\begin{proof}
Assume $H$ is nilpotent with $k\subset Z_1\subset  \cdots\subset Z_t=H$ as in \eqref{nil}.
Let $\gL_{Z_i}$ denote the integral of $Z_i.$ By Lemma \ref{hlamdan}, $\pi_i$ is an algebra isomorphism between $H\gL_{Z_i}$ and  $H_i.$  Since by definition, $\pi_{i}(Z_{i+1})\subset\tl{Z}(H_i)\subset Z(H_i),$ it follows   that $Z_{i+1}\gL_{Z_{i}}$ is central in $H\gL_{Z_{i}}.$

Conversely, assume  $k=N_0\subset N_1\subset\cdots\subset N_t=H$ satisfies the condition in Definition \ref{nilnew}. Let $k\subseteq Z_1\subseteq  \cdots$ be an ascending central series for $H.$

By assumption  we have $N_1\subset Z(H).$ Since $N_1$ is a normal left coideal subagebra of $H,$ it follows from e.g \cite[Prop.1.1]{cwtransactions}, that $N_1\subset \tl{Z}(H)=Z_1.$
Assume by induction  that $N_i\subseteq Z_i.$ We show that $N_{i+1}\subseteq Z_{i+1}.$ Indeed, by definition of $Z_i,$ we have  $\pi_{i}(N_i)\subseteq\pi_{i}(Z_i)=k.$ Since by assumption $N_{i+1}\gL_{N_i}\subset Z(H\gL_{N_i})$ and since $\pi_{i}(\gL_{N_i})=1,$ it follows that
$$\pi_{i}(N_{i+1})=\pi_{i}(N_{i+1}\gL_{N_i})\subset \pi_i(Z(H\gL_{N_i}))=Z(\pi_{i}(H))=Z(H_{i}).$$
But $\pi_{i}(N_{i+1})$ is a normal left coideal subalgebra of $H_{i}$ contained in its center, hence  $\pi_{i}(N_{i+1})\subset \tl{Z}(H_{i}).$ This implies in turn that $\pi_{i+1}(N_{i+1})=k$ and thus $N_{i+1}\subseteq Z_{i+1}.$
Since $N_t=H$ we obtain $Z_t=H$ and thus $H$ is nilpotent.
\end{proof}
As a direct corollary we obtain,
\begin{coro}\label{gen}
Semisimple nilpotent Hopf algebras are solvable.
\end{coro}
\begin{proof}
Assume $H$ has a series as in Definition \ref{nilnew}. Since each $N_i$ is normal in $H,$ it follows that $\gL_{N_i}\in Z(H),$ hence satisfies (i) in Definition \ref{ds} of solvability. Now, since $\gL_{N_i}$ is a central idemoptent of $H,$ we have by assumption $\sum a_1b(Sa_2)\gL_{N_i}=\lla\gep,a\rla b\gL_{N_i}$ for all $a\in H,\,b\in N_{i+1}.$ In particular (ii) of solvability follows.
\end{proof}

Another generalization from groups to Hopf algebras is the following.
\begin{proposition}\label{HK}
Let $K$ be a normal Hopf subalgebra of $H.$ Then $H$ is solvable if and only if $K$ and $H//K$ are solvable.
\end{proposition}
\begin{proof} If $K$ and $H//K$ are solvable then $H$ is solvable  where its solvable series is the union of the solvable series for $K$ and the solvable series in $H$ starting with $K$ as described in Proposition \ref{HN}.

Assume now $H$ is solvable. Clearly $H//K$ is solvable. We want to show that $K$ is solvable. Let $$k=N_0\subset N_1\subset\dots\subset N_t=H$$ be a solvable series for $H.$ We wish to show that the series
$$k=N_0\subset (N_1\cap K)\dots\subset (N_{t}\cap K)= K$$
is a solvable series for $K$ (maybe with repetitions). By Lemma \ref{integrals} $\gl_B\gl_{B_i}$ is an integral for $BB_i$ where $B=(H^*)^K$ and $B_i=(H^*)^{N_i}.$ By \eqref{cap} and Lemma \ref{basic} we have that
$$\gL_{K\cap N_i}=\gl_B\gl_i\rightharpoonup\gL=\gl_B\rightharpoonup\gL_{N_i}.$$
Now, if $a\in N_{i+1}\cap K,$ then we have
$$a\gL_{K\cap N_i}=a(\gl_B\rightharpoonup\gL_{N_i})=\gl_B\rightharpoonup (a\gL_{N_i})=\gl_B\rightharpoonup (\gL_{N_i}a)=(\gl_B\rightharpoonup \gL_{N_i})a=\gL_{K\cap N_i}a.$$
The second and the fourth equality both follow since $B$ is a Hopf subalgebra acting trivially on $a\in K.$
This proves property (i) of solvability.

As for (ii), let $a,b\in K\cap N_{i+1},$ then
\begin{eqnarray*}\lefteqn{(a\ad b)(\gL_{K\cap N_i})=}\\
&=&(a\ad b)(\gl_B\rightharpoonup \gL_{N_i})=\gl_B\rightharpoonup ((a\ad b)\gL_{N_i})=\gl_B\rightharpoonup (\langle\gep,a\rangle b\gL_{N_i})=\langle\gep,a\rangle b\gL_{K\cap N_i}.\end{eqnarray*}
As before, the third and the last equality follow since  $B$ acts trivially on $a\ad b\in K$ and on $b\in K.$
\end{proof}

\bigskip

We are ready to prove the Hopf algebra analogue of Burnside's $p^aq^b$ theorem.

\begin{theorem}\label{burnside}
 Let $H$ be a quasitriangular semisimple  Hopf algebra of dimension $p^aq^b$ over a field $k$ of characteristic $0.$  Then   $H$ is solvable. Moreover, if $N$ is a left coideal subalgebra of $H$ then $H$ has a solvable series containing $N.$
\end{theorem}
\begin{proof}
The proof is by induction on $\dim H.$ Clearly if $H$ is quasitriangular then so are all its homomorphic Hopf images. It is known that the character algebras of quasitriangular Hopf algebras are commutative.

If $\dim H=p$ then $H=k\Z_p,$ a commutative Hopf algebra hence solvable.

Let $H$ be of dimension $p^aq^b,\;a+b>1.$ If $H$ is commutative, we are done. Otherwise, by \cite[Th.4.9]{cwjpaa}, $H$ has a proper non-trivial normal left coideal subalgebra. Let $k\ne N\subset H$ be a minimal normal left coideal subalgebra of $H.$


Consider $N$ as a $D(H)$-module, where $D(H)$ is the Drinfeld double of $H.$ Let $M={\rm LKer}_{_{D(H)}N}\subset D(H).$ Note that since $k\ne (H^*)^N\subset H^*,$ it follows that $k\ne (H^*)^N\tie k\subset M.$
Thus $M$ is a non-trivial normal left coideal subalgebra of $D(H).$ It follows from \cite[Th.4.10]{cwjpaa} that $M$ contains a central grouplike element $\tau\tie\gs\ne \gep\tie 1.$ That is,
\begin{equation}\label{dh}(\gs\ad n)\leftharpoonup\tau\minus=n\end{equation}
for all $n\in N.$
\begin{itemize}
\item If $\tau$ does not act trivially on $N$ then By \eqref{tauh'},  $N\not\subset H',$  the commutator subalgebra of $H.$  Minimality of $N$ implies that $k=N\cap H'.$ But $H//H'$ is a commutative Hopf algebra, hence
    $\pi(h\ad n)=\pi(h)\ad\pi(n) =\langle\gep,h\rangle \pi(n)$ for all $h\in H.$ By Proposition \ref{embedding}, $\pi_{|N}$ is injective, hence $H$ acts-$\ad$ trivially on $N,$ implying that $N\subset Z(H).$ In particular $N$ acts-$\ad$ trivially on itself.

By induction $H//N$ is solvable. Proposition \ref{HN} and the fact that $N$ acts-$\ad$ trivially on itself imply that $H$ is solvable with $k\subset N\subset \dots\subset H.$

\item If $\tau$  acts trivially on $N$ then by \eqref{dh}, $\gs\ad n=n$ for all $n\in N.$ We then have:
\ben
\item If $\gs\ne 1$  then  $k\ne {\rm LKer}_{(N,\ad)}\subset H.$ Minimality of $N$  implies that either $N\subset {\rm LKer}_{(N,\ad)}$ or $N\cap {\rm LKer}_{(N,\ad)}=k.$

If $N\subset {\rm LKer}_{(N,\ad)}$ then it acts-$\ad$ trivially on itself.  By induction $H//N$ is solvable and again Proposition \ref{HN} and the fact that $N$ acts-$\ad$ trivially on itself imply that $H$ is solvable with $k\subset N\subset \dots\subset H.$

 Assume   $N\cap {\rm LKer}_{(N,\ad)}=k.$ By induction $\pi(H)=H//{\rm LKer}_{(N,\ad)}$ is solvable where $\pi(N)$ appears in the series.  Thus we have
$$k\subset \pi (N_1)\subset\dots\subset \pi(N)\subset\dots\subset \pi(H)$$ is a solvable series in $H//{\rm LKer}_{(N,\ad)}.$ By Proposition \ref{embedding}  $\pi_{|N}$ is injective, hence we have that $$k\subset N_1\subset\dots\subset N_l=N$$  is a solvable series in $H.$ Since $H//N$ is solvable by induction, we have by above and Proposition \ref{HN} that $H$ is solvable with
$$k\subset N_1\subset\dots\subset N\subset\dots\subset H.$$
\item If $\gs=1$ then $\tau\tie 1$ is a central grouplike element in $D(H),$ in particular $\gep\ne\tau\in Z(H^*).$
    Since $(H,R)$ is quasitriangular there is a Hopf algebra projection  $\pi_R:D(H)\rightarrow H$  given by
    $$\pi_R(p\tie h)=f_R(p)h$$
     where $R=\sum R^1\ot R^2$ and $f_R(p)=\sum\langle p,R^1\rangle R^2.$

     If $f_R(\tau)\ne 1$ then $f_R(\tau)$ is a central grouplike element of $H,$ hence the Hopf center of $H,\,\tl{Z}(H)\ne k.$ By induction $H//\tl{Z}(H)$ is solvable, hence by Proposition \ref{HK}, $H$ is solvable with the series
     $$k\subset \tl{Z}(H)\subset\dots\subset H.$$

    By the same argument, if $f_{R^t}(\tau)=\sum\langle\tau,R^2\rangle R^1\ne 1$ then $H$ is solvable.

    Assume now $f_R(\tau)=f_{R^t}(\tau)=1.$   Clearly $N\subset K=H^{k\langle\tau\rangle}$ since by assumption $\tau$ acts trivially on $N.$ Moreover, $K$ is a normal Hopf subalgebra of $H$ since $k\lla\tau\rla$ is a normal Hopf subalgebra of $H^*.$ Once we show that $K$ is quasitriangular, we are done by induction and Proposition \ref{HK}.

    Let $K^*=H^*//k\langle\tau\rangle=\ol{H^*}.$ We have a co-quasitriangular  structure on $K^*$ given by:
    $$\qquad\qquad\langle \ol{p}\,|\,\ol{q}\rangle=\sum \langle p,R^2\rangle\langle q, R^1\rangle=\langle f_R(p),q\rangle=\langle f_{R^t}(q),p\rangle.$$
    Clearly $\langle\,|\,\rangle$ satisfies properties of co-quasitriangularity. We need to show that it is well defined. Indeed, if $\ol{p}=0$ or respectively $\ol{q}=0,$ then $p\in H(k\langle\tau\rangle^+)$ or respectively $q\in H(k\langle\tau\rangle)^+$). Since $f_R$ and $f_{R^t}$ are multiplicative (and respectively anti multiplicative) it follows that $f_R(p)=0$ or respectively $f_{R^t}(q)=0.$ This proves that $\langle\,|\,\rangle$ is well defined, and thus $K$ is quasitriangular.
    This concludes the proof of the theorem.
\een
\end{itemize}
 \end{proof}

\end{document}